\documentclass[draft]{amsart}
\usepackage{amssymb,graphicx}
\usepackage{multirow,comment,caption}

\newcommand{\df}{\dfrac}
\newcommand{\tf}{\tfrac}

\newcommand{\z}{\zeta}

\newcommand{\beqs}{\begin{equation*}}
\newcommand{\eeqs}{\end{equation*}}
\numberwithin{equation}{section}
 \theoremstyle{plain}
\newtheorem{theorem}{Theorem}[section]
\newtheorem{lemma}[theorem]{Lemma}
\newtheorem{corollary}[theorem]{Corollary}

\theoremstyle{remark}

\begin{document}

\makeatletter
\def\imod#1{\allowbreak\mkern10mu({\operator@font mod}\,\,#1)}
\makeatother

\author{Alexander Berkovich}
   \address{Department of Mathematics, University of Florida, 358 Little Hall, Gainesville FL 32611, USA}
   \email{alexb@ufl.edu}

\author{Frank Patane}
   \address{Department of Mathematics, University of Florida, 358 Little Hall, Gainesville FL 32611, USA}
   \email{frankpatane@ufl.edu}

\title[\scalebox{.9}{Essentially Unique Representations by Certain Ternary Quadratic Forms}]{Essentially Unique Representations by Certain Ternary Quadratic Forms}
     
\begin{abstract} 

In this paper we generalize the idea of ``essentially unique'' representations by ternary quadratic forms. We employ the Siegel formula, along with the complete classification of imaginary quadratic fields of class number less than or equal to 8, to deduce the set of integers which are represented in essentially one way by a given form which is alone in its genus. We consider a variety of forms which illustrate how this method applies to any of the 794 ternary quadratic forms which are alone in their genus. As a consequence, we resolve some conjectures of Kaplansky regarding unique representation by the forms $x^2 +y^2 +3z^2$, $x^2 +3y^2 +3z^2$, and $x^2 +2y^2 +3z^2$ \cite{kap}.

\end{abstract}

\keywords{representations of integers, sum of three squares, local densities, ternary quadratic forms, Siegel formula}

 \subjclass[2010]{11B65, 11E16, 11E20, 11E25, 11E41, 11F37}

\date{\today}
   
\maketitle

   \section{Introduction}
\label{intro}

The concept of connecting the number of representations of a ternary quadratic form to the class number for binary quadratic forms dates back to Gauss \cite{gauss}. Gauss was the first to introduce many fundamental concepts such as discriminant, positive definite form, and equivalence of forms. After introducing these fundamental notions, he related the number of representations of an integer by $x^2 + y^2 +z^2$ to what is essentially the class number for binary quadratic forms.\\

Representation of integers by $x^2 +y^2 +z^2$, has been studied by many mathematicians since the time of Gauss. Building on the work of Hardy, Bateman \cite{bate} derived and
proved the formula for the number of representations of a positive integer as the sum of three squares. We point out that this representation formula is a special case of the more general Siegel formula which can be found in \cite{siegel}. We mention this since our treatment often relies on the Siegel formula, which we will describe in the next section.\\

Rather than discussing the total number of representations by a quadratic form, one can identify solutions according to a given relation. In the case of the form $x^2 +y^2 +z^2$, identifying solutions which are the same up to order and sign is equivalent to partitioning a number into three squares. In 1948, Lehmer considered partitions of an integer into $k$ squares \cite{lehm}. We refer the reader to \cite{sel} for a recent (2004) discussion of this topic.\\

In 1984, Bateman and Grosswald essentially classified all integers which have one representation up to order and sign by $x^2 +y^2 +z^2$ \cite{batgros}. Their proof assumed they had the complete list of discriminants of binary quadratic forms with class number less than or equal to 4. In 1992, Arno completely classified all discriminants of binary quadratic forms with class number less than or equal to 4 \cite{arno}. Bateman and Grosswald's assumption was proven correct.\\

In 1997, Kaplansky considered the forms $x^2 +y^2 +2z^2$, $x^2 +2y^2 +2z^2$, and $x^2 +2y^2 +4z^2$ \cite{kap}. He identified solutions which are the same up to ``order and sign'', and deduced which numbers are represented in essentially one way by the aforementioned forms. He utilized the completed list of discriminants of binary quadratic forms with class number less than or equal to 4 to deduce the integers with essentially unique representation by the forms he considered. Kaplansky then conjectured about the numbers which are represented in essentially one way by the forms $x^2 +2y^2 +3z^2$, $x^2 +3y^2 +3z^2$, and $x^2 +y^2 +3z^2$.\\

In this paper we employ the results of Watkins \cite{watkins} to resolve Kaplansky's conjectures. Furthermore, we extend the idea of essentially unique representation beyond diagonal forms, where we must consider more than ``order and sign''. Our treatment will apply to any of the 794 ternary quadratic forms which are alone in their genus. The determination of these forms was first explored by Watson \cite{wat1}, with the final touch delivered by Jagy, Kaplansky, and Schiemann \cite{jagykap}. (See also \cite{lorch}.) For the remainder of this paper, we call a form idoneal when it is alone in its genus. We point out that $x^2 +y^2 +z^2$, along with all the ternary forms discussed in \cite{kap}, are idoneal.\\

In Section \ref{notprim}, we will give the necessary definitions and notation as well as discuss automorphs and ``essentially unique'' representations. We then outline the general approach of how to use the Siegel formula along with class number bounds to derive integers which are represented in essentially one way by an idoneal ternary form. The largest class number bound we utilize in this paper is 8. We have compiled tables in the Appendix which list all discriminants of binary quadratic forms with class number less than or equal to 8 according to class group type. \\

In Section \ref{ex1}, we will consider the non-diagonal forms $x^2 +y^2 +z^2 +yz +xz +xy$,  $3x^2 +3y^2 +3z^2 -2yz +2xz +2xy,$ and $x^2 +3y^2 +3z^2 +2yz$. These forms are selected and grouped together in Section \ref{ex1} because we treat these forms by relating them to $x^2+y^2+z^2$. This generalizes the approach of Kaplansky \cite{kap} to non-diagonal forms, and we comment that the three selected forms are among many which can be handled in a similar fashion. In particular, if $f$ is an idoneal form of discriminant $\Delta=2^k$, then one can find the integers which are uniquely represented by $f$ by reducing $f$ to $x^2+y^2+z^2$.\\

In Section \ref{ex2}, we will examine the non-diagonal forms $5x^2 +13y^2 +20z^2 -12yz +4xz +2xy$ and $7x^2 +15y^2 +23z^2 +10yz +2xz +6xy$. Both of these forms can be treated by the methods of Section \ref{ex1}, however we chose to use the Siegel formula along with local density considerations to derive the integers which they represent in essentially one way. \\
 
In Section \ref{ex3}, we resolve the aforementioned conjectures of Kaplansky. Explicitly, we find the integers which are represented in essentially one way by $x^2 +y^2 +3z^2$, $x^2 +3y^2 +3z^2$, and $x^2 +2y^2 +3z^2$, by applying the method outlined in Section \ref{notprim}.\\ 

Section \ref{concl}, also called the Outlook, contains our concluding remarks. In the Outlook we consider the form $x^2 +3y^2 +3z^2 +yz +xy$ which is not idoneal. We sketch the proof that we have found all integers which are represented in essentially one way by this form. We conclude this paper with prospects for future work.

\section{Notation and Preliminaries}
\label{notprim}

We use the notation $(a,b,c,d,e,f)$ to represent the positive ternary quadratic form $ax^2 +by^2 +cz^2 +dyz +exz +fxy$. We remark that this paper only considers positive ternary quadratic forms. We use $(a,b,c,d,e,f;n)$ to denote the total number of representations of $n$ by $(a,b,c,d,e,f)$. We take $(a,b,c,d,e,f;n)=0$ when $n \not\in \mathbb{N}$. 
The associated theta series to the form $(a,b,c,d,e,f)$ is
\begin{equation}
\label{thet}
\vartheta(a,b,c,d,e,f,q):=\sum_{x,y,z}q^{ax^2 +by^2 +cz^2 +dyz +exz +fxy} = \sum_{n\geq 0}(a,b,c,d,e,f;n)q^n.
\end{equation}
The discriminant $\Delta$ of $(a,b,c,d,e,f)$ is defined as
\[
\Delta:= \frac{1}{2} \mbox{det}\left( \begin{array}{lll}
 2a& f &e  \\
f & 2b &d  \\
 e&d  &2c  \\\end{array} \right) = 4abc +def -ad^2 -be^2 -cf^2.
\]
We note that the discriminant $\Delta>0$ for a positive ternary quadratic form.
Two ternary quadratic forms of discriminant $\Delta$ are in the same genus if they are equivalent over $\mathbb{Q}$ via a transformation matrix in $SL(3,\mathbb{Q})$ whose entries have denominators coprime to $2\Delta$.\\

	 Let $A$ be a 3 by 3 matrix of determinant $\pm 1$. $A$ is an automorph for the form $(a,b,c,d,e,f)$ if the action of $A$ on $(a,b,c,d,e,f)$ leaves $(a,b,c,d,e,f)$ unchanged. We denote the set of automorphs of $(a,b,c,d,e,f)$ by Aut$(a,b,c,d,e,f)$. A discussion of automorphs for ternary quadratic forms is given in \cite{dickson}. We use Sage 5.1 to explicitly compute the automorphs for the forms considered in this paper. We now give a brief example by considering the 8 automorphs of the form $(1,3,4,3,1,0)$. We have	
	\begin{align*}
	\mbox{Aut}(1,3,4,3,1,0)= 
\Bigg\{&\left( \begin{array}{lll}
 $1$& $0$ &$0$  \\
$0$ & $1$ &$0$  \\
 $0$&$0$  &$1$  \\\end{array} \right),\left( \begin{array}{lll}
 $-$1$$& $0$ &$0$  \\
$0$ & $-$1$$ &$0$  \\
 $0$&$0$  &$-$1$$  \\\end{array} \right),\left( \begin{array}{lll}
 $1$& $0$ &$1$  \\
$0$ & $-$1$$ &$0$  \\
 $0$&$0$  &$-$1$$  \\\end{array} \right),\left( \begin{array}{lll}
 $-$1$$& $0$ &$-$1$$  \\
$0$ & $1$ &$0$  \\
 $0$&$0$  &$1$  \\\end{array} \right),\\
& \left( \begin{array}{lll}
 $1$& $0$ &$0$  \\
$0$ & $-$1$$ &$-$1$$  \\
 $0$&$0$  &$1$  \\\end{array} \right),\left( \begin{array}{lll}
 $-$1$$& $0$ &$0$  \\
$0$ & $1$ &$1$  \\
 $0$&$0$  &$-$1$$  \\\end{array} \right),\left( \begin{array}{lll}
 $1$& $0$ &$1$  \\
$0$ & $1$ &$1$  \\
 $0$&$0$  &$-$1$$  \\\end{array} \right),\left( \begin{array}{lll}
 $-$1$$& $0$ &$-$1$$  \\
$0$ & $-$1$$ &$-$1$$  \\
 $0$&$0$  &$1$  \\\end{array} \right)\Bigg\}.
\end{align*}
	
	To give a further illustration, we note that $(1,3,4,3,1,0;19)=12$. Under the action of \\Aut$(1,3,4,3,1,0)$, the solutions form two orbits:
	\begin{align*}
	O_1:&=\{(-4, -1, 0), (-4, 1, 0), (4, -1, 0), (4, 1, 0)\},\\
	O_2:&=\{(-3, -2, 2), (-3, 0, 2), (-1, 0, -2), (-1, 2, -2), (1, -2, 2), (1, 0,2), (3, 0, -2), (3, 2, -2)\}.
	\end{align*}
The solutions in $O_1$ are easily identified as the solution $(4,1,0)$ up to sign. However the solutions in $O_2$ are not so readily identified as being equivalent under the action of automorphs.\\

Identifying solutions which are equivalent under the action of automorphs is the way to generalize previous authors' (\cite{batgros}, \cite{kap}, \cite{lehm}) notion of solutions being equivalent up to ``order and sign''.\\
When the solutions form exactly $k$ orbits under the action of automorphs, we say the form represents the integer in essentially $k$ ways. We say an integer has an essentially unique representation when the solutions form 1 orbit under the action of automorphs. We also note that if $f$ is any ternary quadratic form then $f(x,y,z)=n$ implies $f(-x,-y,-z)=n$. Thus if $f$ represents a positive integer $n$, then $f$ represents $n$ in at least two ways. Hence there is little ambiguity if we say that $f$ uniquely represents an integer when the solutions form 1 orbit under the action of automorphs.\\

 The focus of this paper is concerned with finding integers which have an essentially unique representation by a given idoneal form. In particular, we give a method which enables one to find all integers which are represented in essentially one way, by an idoneal ternary quadratic form. We now introduce an essential tool to our method, the celebrated Siegel theorem for positive ternary quadratic forms.

\begin{theorem}
\label{sw}
Let $G$ be a genus of positive ternary quadratic forms of discriminant $\Delta$. Then 
\begin{equation}
\sum_{t \in G} \frac{R_t(n)}{|\textnormal{Aut}(t)|} = 4\pi M(G) \sqrt{\frac{n}{\Delta}}\prod_{p} d_{G,p}(n),
\end{equation}
\end{theorem}
\noindent
where $R_t(n)$ denotes the total number of representations of $n$ by $t$, and it is understood that the sum on the left is over representatives of each equivalence class in the genus $G$. The product on the right is over all primes $p$, and the mass of $G$ is defined as
\[
M(G):=\sum_{t \in G} \frac{1}{|\mbox{Aut}(t)|}.
\]
Let $t:=ax^2 +by^2 +cz^2 +dyz +exz +fxy$ be any form in $G$. Then the $p$-adic local density $d_{G,p}$ (also called $d_{t,p}$) is
\begin{equation}
\label{locald}
d_{t,p}(n):=\lim_{k\to \infty} p^{-2k}|\{(x,y,z)\in \mathbb{Z}^3 : ax^2 +by^2 +cz^2 +dyz +exz +fxy\equiv n \imod{p^k}\}|.
\end{equation}
We remark that the limit in \eqref{locald} can be removed as long as we take $k$ large.\\

	Theorem \ref{sw} is a special case of the general Siegel theorem given in \cite{siegel}. When the form $f$ is idoneal, Theorem \ref{sw} gives an explicit formula for $R_f(n)$.
	\begin{corollary}
\label{sws}
Let $f$ be an idoneal ternary quadratic form of discriminant $\Delta$. Then 
\[
R_{f}(n) = 4\pi \sqrt{\df{n}{\Delta}} \prod_{p} d_{f,p}(n),
\]
with all notation as previous.
\end{corollary}
\noindent
In \cite{siegel} Siegel shows that when $(2\Delta,p)=1$ we have
\begin{equation}
\label{exsi}
d_{t,p}(n) = \left\{ \begin{array}{ll} 1+\frac{1}{p} + \frac{1}{p^{k+1}}\left( \left(\frac{-m\Delta}{p}\right)-1\right)&  n=mp^{2k}, p \nmid m,\\
&\\
				\left(\frac{1}{p}+1\right)\left(1-\frac{1}{p^{k+1}}\right)&  n=mp^{2k+1}, p \nmid m,\\
     \end{array}
     \right.
	\end{equation}
	where $\left(\tf{r}{p}\right)$ is the Legendre symbol.
	When we use Corollary \ref{sws} along with equation \eqref{exsi}, we obtain a very explicit formula for the total number of representations of an integer by an idoneal form. Employing \eqref{exsi} it can be shown that
	\begin{equation}
	\label{sp}
	\prod_{p\nmid 2\Delta} d_{f,p}(n) = \df{8}{\pi^2}L(1,\chi(\Delta n)) P(n,\Delta) \prod_{2<p \mid \Delta}\df{1}{1-\tf{1}{p^2}},
	\end{equation}
where $L(1,\chi(n))$ is given by
	\begin{equation}
	\label{el}
	L(1,\chi(n)):=\sum_{m=1}^{\infty}\df{\left(\tf{-4n}{m}\right)}{m}=\prod_{p>2}\df{1}{1-\tf{\left(\tf{-n}{p}\right)}{p}},
	\end{equation}
	and $\chi(n):=\left(\frac{-4n}{\bullet}\right)$. Lastly $P(n,\Delta)$ is given by the finite product
	\begin{equation}
	\label{p}
P(n,\Delta) := \prod_{\substack{(p^{2})^b\mid\mid n,\\p \nmid 2\Delta}}\left( 1+ \df{1}{p} + \df{1}{p^2}+\cdots +\df{1}{p^{b-1}} + \df{1}{p^b (1 - \left(\tf{-\Delta np^{-2b}}{p}\right)p^{-1})}\right),
	\end{equation}
  where the product is over all primes $p \nmid 2\Delta$ such that $p^2\mid n$, and $b$ is the largest integer such that $p^{2b}\mid n$. We note that the only property of $P(n,\Delta)$ that we use is $P(n,\Delta)$ is a finite product with $1\leq P(n,\Delta)$. Lastly, $P(n,4^k)=P(n,4)=P(n,1)$ and so we define the abbreviated $P(n):=P(n,1)$.\\
	
	Combining Corollary \ref{sws} with \eqref{sp} yields
	\begin{equation}
	\label{mf}
	R_{f}(n) = \frac{32\sqrt{n}}{\pi\sqrt{\Delta}} L(1,\chi(\Delta n))\cdot P(n,\Delta)\cdot \prod_{p\mid 2\Delta} d_{f,p}(n)\prod_{2<p \mid \Delta}\df{1}{1-\tf{1}{p^2}},
	\end{equation}
	where $f$ is an idoneal ternary quadratic form of discriminant $\Delta$.\\
	
	If we factor $n$ as $n=4^a \cdot m\cdot  d^2$ with $2 \nmid d$ and $m$ squarefree, then
	\begin{align}
\label{lphi}
L(1,\chi(n)) &=L(1,\chi(m))\prod_{p\mid d}1-\frac{\left(\tf{-m}{p}\right)}{p}.
\end{align}

	Dirichlet gives a wonderful connection between $L(1,\chi(m))$ and $h(D)$, the number of reduced primitive binary quadratic forms of discriminant $D=-m$ or $D=-4m$ \cite{dir}. We call $h(D)$ the class number of discriminant $D$. The relationship between $L(1,\chi(m))$ and $h(D)$ is given in the following theorem.
	\begin{theorem}
	\label{tthg}
	 For $m$ squarefree and $\chi(m):=\left(\frac{-4m}{\bullet}\right)$, we have
\begin{equation}
\label{h}
L(1,\chi(m)) = \left\{ \begin{array}{ll}
        \dfrac{\pi}{4}&  m=1,\\
				\vspace{-.2cm}&\\
				\dfrac{\pi}{2\sqrt{3}} &  m=3,\\
				\vspace{-.2cm}&\\
				\dfrac{3\pi}{2\sqrt{m}}h(-m) &  3<m\equiv 3 \imod{8},\\
				\vspace{-.2cm}&\\
				\dfrac{\pi}{2\sqrt{m}}h(-m) &  m \equiv 7 \imod{8},\\
				\vspace{-.2cm}&\\
				\dfrac{\pi}{2\sqrt{m}}h(-4m) & 1< m\equiv 1,2 \imod{4}.\\
     \end{array}
     \right.
\end{equation}
		\end{theorem}
		\noindent
	See \cite{buell} for details.\\
	
	 Let $f$ be an idoneal ternary quadratic form. A necessary but not sufficient condition for $f$ to uniquely represent $n$ up to the action of automorphs, is
\begin{equation}
\label{formula}
0<R_{f}(n) \leq |\text{Aut}(f)|.
\end{equation}
 We employ \eqref{mf} along with Theorem \ref{tthg} to give an explicit lower bound for $R_f(n)$ in terms of the class number. We then classify the $n$ which satisfy \eqref{formula}, and call this set the prelist of $f$, denoted by $\text{Prelist}(f)$. We need only check which elements of $\text{Prelist}(f)$ have the property that their solutions form one orbit under Aut$(f)$. There are only a finite number of elements of $\text{Prelist}(f)$ we need to check, and we employ Maple V.15 to compute the number of orbits of the solutions. Of course, solving for the elements of $\text{Prelist}(f)$ requires one to have information on the bounds of the class number. Hence we now discuss class number considerations.\\

Attempting to solve $h(d)=k$ dates back to Gauss \cite[Section V, Article 303]{gauss}, where Gauss conjectured the list of imaginary quadratic fields of small class number. 
The case $h(d)=1$ was essentially first completed by Heegner in 1951 \cite{heg}. In 1967, Baker and Stark gave independent proofs for the $h(d)=1$ case as well. Many mathematicians have done extensive work towards solving $h(d)=k$ for $k\leq 8$. (See \cite{arno}, \cite{baker}, \cite{gold}, and \cite{gross}.)\\
	
	The latest developments are given by Watkins \cite{watkins}. According to Watkins, the largest (in magnitude) fundamental discriminant $d$ with class number less than or equal to $8$, is $d=-6307$.  We can enumerate those with non-fundamental discriminant by employing the formula
	\begin{equation}
	\label{veq}
	h(D) = h(d) \cdot f\cdot  \frac{w_D}{w_d}  \prod_{p \mid f} 1 - \tf{\left(\tf{d}{p}\right)}{p},
	\end{equation}
	where $D=d\cdot f^2$, $f$ is the conductor of $D$, and $w_\mathfrak{D}$ is given by
	\[
	w_\mathfrak{D} :=\left\{ \begin{array}{ll}
        6&  \mathfrak{D}  =-3,\\
				   4&  \mathfrak{D}  =-4,\\
					   2&  \mathfrak{D}  <-4.\\
     \end{array}
     \right.
	\]
	Equation \eqref{veq} is Lemma 2.13 in \cite{voight}. An application of \eqref{veq} is $h(d\cdot f^2)\leq 8$ and $d=-3$, $d=-4$, $|d| > 4$, implies $f \leq 90, 60, 30$, respectively.\\

Another important use of \eqref{veq}, is we can combine \eqref{lphi} with \eqref{veq} to remove the restriction of $m$ being squarefree in Theorem \ref{tthg}.\\

In the Appendix, we include tables of the 527 discriminants (fundamental and non-fundamental) with class number $\leq 8$ organized by isomorphism class of the class group. We generated this complete set by utilizing the bounds found in \cite{watkins} along with \eqref{veq}. We then used PARI/GP V.2.7.0 to identify the isomorphism class of each class group, and compile the tables.\\

We now move on to Section \ref{ex1}, where we consider the form $(1,1,1,0,0,0)$ and derive the corresponding prelist. We then use this prelist to find the integers which are uniquely represented by the forms $(1,1,1,1,1,1), (3,3,3,-2,2,2),$ and $(1,3,3,2,0,0)$.

\section{Some ternary forms of discriminant $2,4,32,64$ and their relation to $x^2+y^2+z^2$}
\label{ex1}

In \cite{bate}, Bateman shows
\begin{equation}
\label{thm}
(1,1,1,0,0,0;n)=\dfrac{16\sqrt{n}}{\pi}d_{s,2}(n) L(1,\chi(n))P(n),
\end{equation}
where $s:=x^2 +y^2 +z^2$, and all other notation is given in Section \ref{notprim}. We comment that equation \eqref{thm} follows from \eqref{mf} as well.

Since $|$Aut$(s)|=48$, the prelist of $s$ consists of all $n$ with $0<(1,1,1,0,0,0;n)\leq 48$. By congruence considerations it is easy to see that $(1,1,1,0,0,0;n)=(1,1,1,0,0,0;4n)$, and so we restrict to $4 \nmid n$.\\


We refer to \cite[Equation (1.5)]{berk}, for the function $d_{s,2}(n)$. We have
\begin{equation}
\label{psi}
d_{s,2}(n) = \left\{ \begin{array}{ll}
3/2 &  n\equiv 1,2\imod{4}, \\
1 &  n\equiv 3\imod{8}, \\
0 &   n\equiv 7\imod{8}, \\
     \end{array}
     \right.
\end{equation}
 for $n \not\equiv 0 \imod{4}$. We need not consider $n \equiv 7 \imod{8}$ since \eqref{psi} and \eqref{thm} imply\\ $(1,1,1,0,0,0;n)=0$ for such $n$.\\

Both $n=1,3$ satisfy $0<(1,1,1,0,0,0;n)\leq 48$. Employing \eqref{h}, \eqref{thm}, and \eqref{psi}, we find
\begin{equation}
\label{sss1}
(1,1,1,0,0,0;n)\geq \left\{ \begin{array}{ll}
				24h(-n) &  3<n\equiv 3 \imod{8},\\
			12h(-4n) & 1< n\equiv 1,2 \imod{4}.\\
     \end{array}
     \right.
\end{equation}

Our goal is to solve for $n$ which satisfy
\begin{equation}
\label{sss2}
48\geq \left\{ \begin{array}{ll}
				24h(-n) &  3<n\equiv 3 \imod{8},\\
			12h(-4n) & 1< n\equiv 1,2 \imod{4}.\\
     \end{array}
     \right.
\end{equation}
We point out that to solve \eqref{sss2} we need information regarding discriminants with class number 4. Employing the tables in the Appendix, we solve \eqref{sss2} and find there are a total of 53 integers $n$ with $4 \nmid n$, and $0<(1,1,1,0,0,0;n) \leq 48$. We remark that three spurious solutions, $n=49, 75, 99$, satisfy \eqref{sss2} yet $(1,1,1,0,0,0;n) > 48$. $\text{Prelist}(1,1,1,0,0,0)$ consists of the integers $4^k \cdot v$, $k \geq 0$, and\\ $v\in$ 
\{1, 2, 3, 5, 6, 9, 10, 11, 13, 14, 17, 18, 19, 21, 22, 25, 27, 30, 33,
34, 35, 37, 42, 43, 46, 51, 57, 58, 67, 70, 73, 78, 82, 85, 91, 93,
97, 102, 115, 123, 130, 133, 142, 163, 177, 187, 190, 193, 235, 253,
267, 403, 427\}.

In Table \ref{prel} we tabulate the integers $4 \nmid n$ with $0<(1,1,1,0,0,0;n) \leq k$ for select values of $k$.

\begin{table}[htb]  
	\caption{} \label{prel}
     \begin{center} 
       \begin{tabular}{|l|l|}              \hline 
			$(1,1,1,0,0,0;n)$ & $4\nmid n$\\ \hline
			$0<(1,1,1,0,0,0;n)<6$ & $n \in \{\}$\\ \hline
			$(1,1,1,0,0,0;n)=6$ & $n \in \{1\}$\\ \hline
			$6<(1,1,1,0,0,0;n)\leq 8$ & $n \in \{3\}$\\ \hline
			$8<(1,1,1,0,0,0;n)\leq 12$ & $n \in \{2\}$\\ \hline
			$12<(1,1,1,0,0,0;n)< 24$ & $n \in \{\}$\\ \hline
			$(1,1,1,0,0,0;n)= 24$ & $n \in \{5, 6, 10, 11, 13, 19, 22, 37, 43, 58, 67, 163\}$\\ \hline
			$24<(1,1,1,0,0,0;n)< 48$ & $n \in \{9, 18, 25, 27\}$\\ \hline
			$(1,1,1,0,0,0;n)=48$ & $n \in \{$14, 17, 21, 30, 33, 34, 35, 42, 46,\\
			& 51, 57, 70, 73, 78, 82, 85, 91, 93, 97, 102, 115, 123, \\
			& 130, 133, 142, 177, 187, 190, 193, 235, 253, 267, 403, 427$\}$\\ \hline
			\end{tabular}
			\begin{flushright}
			 .
			 \end{flushright}
     \end{center}
   \end{table}

We now check which integers in $\text{Prelist}(1,1,1,0,0,0)$ are represented in an essentially unique way, and arrive at the following theorem.

\begin{theorem}
\label{ttz}
The form $(1,1,1,0,0,0)$ uniquely represents $n$ (up to action of automorphs) if and only if $n=4^k\cdot v$, $k\geq 0$, and \\ $v \in$
\{$1$, $2$, $3$, $5$, $6$, $10$, $11$, $13$, $14$, $19$, $21$, $22$, $30$, $35$, $37$, $42$, $43$, $46$, $58$, $67$, $70$, $78$, $91$, $93$, $115$, $133$, $142$, $163$, $190$, $235$, $253$, $403$, $427$\}.
\end{theorem}

Theorem \ref{ttz} was previously established in \cite{batgros}. However we additionally derived Table \ref{prel} in the process of proving Theorem \ref{ttz}, and we employ Table \ref{prel} to find the integers which are uniquely represented by certain forms that are considered in subsequent sections.\\

The idea of relating a form to $x^2+y^2+z^2$ to count the number of representations is not new. Indeed, this idea was developed by Dickson \cite{mdickson}, and is the main technique of Kaplansky \cite{kap} where he determined the integers which are represented in an essentially unique way by the forms $x^2+y^2+2z^2$, $x^2+2y^2+2z^2$, and $x^2+2y^2+4z^2$. We mention that although all the Theorems listed in \cite{kap} are correct, some of the main Lemmas (Lemma 3.1--3.3) can easily be misinterpreted. We chose to restate and augment these important Lemmas of Kaplansky.

\begin{lemma}[Kaplansky]
\label{kap1l}
Let $n$ be a nonnegative integer. We have
\begin{equation}
\label{kap1}
(1,1,2,0,0,0;2n)= (1,1,1,0,0,0;n),
\end{equation}
\begin{equation}
\label{kap1b}
3(1,1,2,0,0,0;2n+1)= (1,1,1,0,0,0;4n+2).
\end{equation}
\end{lemma}

\begin{lemma}[Kaplansky]
\label{kap2l}
Let $n$ be a nonnegative integer. We have
\begin{equation}
\label{kap2}
(1,2,2,0,0,0;4n)= (1,1,1,0,0,0;n),
\end{equation}
\begin{equation}
\label{kap2b}
3(1,2,2,0,0,0;4n+1)= (1,1,1,0,0,0;4n+1),
\end{equation}
\begin{equation}
\label{kap2c}
3(1,2,2,0,0,0;4n+2)= (1,1,1,0,0,0;4n+2),
\end{equation}
\begin{equation}
\label{kap2d}
(1,2,2,0,0,0;4n+3)= (1,1,1,0,0,0;4n+3).
\end{equation}
\end{lemma}
\begin{lemma}[Kaplansky]
\label{kap3l}
Let $n$ be a nonnegative integer. We have
\begin{equation}
\label{kap3}
(1,2,4,0,0,0;8n)= (1,1,1,0,0,0;n),
\end{equation}
\begin{equation}
\label{kap3c}
3(1,2,4,0,0,0;8n+2)= (1,1,1,0,0,0;4n+1),
\end{equation}
\begin{equation}
\label{kap3d}
3(1,2,4,0,0,0;8n+4)= (1,1,1,0,0,0;4n+2),
\end{equation}
\begin{equation}
\label{kap3e}
(1,2,4,0,0,0;8n+6)= (1,1,1,0,0,0;4n+3),
\end{equation}
\begin{equation}
\label{kap3b}
6(1,2,4,0,0,0;2n+1)= (1,1,1,0,0,0;4n+2).
\end{equation}
\end{lemma}

We now go beyond Kaplansky's forms and relate the non-diagonal forms $(1,1,1,1,1,1),$\\ $(3,3,3,-2,2,2)$, and $(1,3,3,2,0,0)$ to $x^2+y^2+z^2$ to find the integers which they uniquely represent. \\

The form $(1,1,1,1,1,1)$ has discriminant 2 and 48 automorphs. We write $f(x,y,z):=x^2 +y^2+z^2 +yz +xz +xy$. It is easy to check that
\[
f(x,y,z) \equiv 0 \imod{2},
\]
if and only if $x \equiv y \equiv z \imod{2}$.\\

For any $x,y,z$ we must have $x \equiv y \equiv z \imod{2}$ or we have one of the following: $x \equiv y \not\equiv z \imod{2}$, $x \equiv z \not\equiv y \imod{2}$, or $y \equiv z \not\equiv x \imod{2}$. We note that
\begin{equation}
\label{lee}
\sum_{x\equiv y \not\equiv z\imod{2}}q^{f(x,y,z)}=\sum_{x\equiv z \not\equiv y\imod{2}}q^{f(x,y,z)}=\sum_{y\equiv z \not\equiv x\imod{2}}q^{f(x,y,z)},
\end{equation}
and thus we have
\begin{equation}
\label{dis1}
\sum_{x,y,z}q^{f(x,y,z)}=\sum_{x\equiv y \equiv z\imod{2}}q^{f(x,y,z)} + 3\sum_{x\equiv y \not\equiv z\imod{2}}q^{f(x,y,z)}.
\end{equation}

The substitution $x \mapsto (-x +y+z)$, $y \mapsto (x-y+z)$, and $z \mapsto (x+y-z)$, guarantees the condition $x \equiv y \equiv z \imod{2}$, hence
\begin{equation}
\label{dis2}
\sum_{x\equiv y \equiv z\imod{2}}q^{f(x,y,z)}=\sum_{x,y,z}q^{f(-x +y+z,x-y+z,x+y-z)}=\sum_{x,y,z}q^{2(x^2+y^2+z^2)}.
\end{equation}

The substitution $x \mapsto (-x +y+z)$, $y \mapsto (x-y+z)$, and $z \mapsto (x+y-z+1)$, gives $x\equiv y \not\equiv z\imod{2}$, and we have
\begin{equation}
\label{dis22}
\sum_{x\equiv y \not\equiv z\imod{2}}q^{f(x,y,z)}=\sum_{x,y,z}q^{f(-x +y+z,x-y+z,x+y-z+1)}=q\sum_{x,y,z}q^{2(x^2+y^2+z^2)+2(y+z)}.
\end{equation}

We have now proven the following lemma.
\begin{lemma}
\label{tot}
\begin{equation}
\sum_{x,y,z}q^{x^2 +y^2+z^2 +yz +xz +xy} = \sum_{x,y,z}q^{2(x^2+y^2+z^2)} + 3q\sum_{x,y,z}q^{2(x^2+y^2+z^2)+2(y+z)}.
\end{equation}
\end{lemma}

Lemma \ref{tot} implies 
\begin{equation}
\label{odd1}
(1,1,1,1,1,1;2n+1) = (1,1,1,0,0,0;4n+2),
\end{equation}
and
\begin{equation}
\label{even1}
(1,1,1,1,1,1;2n) = (1,1,1,0,0,0;n),
\end{equation}
for any nonnegative integer $n$.\\
Equations \eqref{odd1} and \eqref{even1} relate $(1,1,1,1,1,1;n)$ to $(1,1,1,0,0,0;n)$ for any nonnegative integer. We can use Table \ref{prel} to determine the solutions to 
\begin{equation}
\label{oddl}
0<(1,1,1,1,1,1;2n+1)=(1,1,1,0,0,0;4n+2) \leq 48,
\end{equation}
and
\begin{equation}
\label{evenl}
0<(1,1,1,1,1,1;2n)=(1,1,1,0,0,0;n) \leq 48.
\end{equation}

We find the solutions to $0<(1,1,1,1,1,1;2n) \leq 48$ to be $2n =4^k \cdot v$, $k\geq 0$, and $v$ is in the set
\{2, 4, 6, 10, 12, 18, 20, 22, 26, 28, 34, 36, 38, 42, 44, 50, 54, 60, 66, 68, 70, 74, 84, 86, 92, 102, 114, 116, 134, 140, 146, 156, 164, 170, 182, 186, 194, 204, 230, 246, 260, 266, 284, 326, 354, 374, 380, 386, 470, 506, 534, 806, 854\}.

We find the solutions to $0<(1,1,1,1,1,1;2n+1) \leq 48$ to be 
\begin{equation}
\label{oddsol}
2n+1 =1, 3, 5, 7, 9, 11, 15, 17, 21, 23, 29, 35, 39, 41, 51, 65, 71, 95.
\end{equation}
We use Maple V.15 to check the above candidates for unique representation, and arrive at the following theorem.
\begin{theorem}
\label{333oe}
The form $(1,1,1,1,1,1)$ uniquely represents the integer $n$ (up to action of automorphs) if and only if  $n=4^k\cdot v$, $k\geq 0$, and \\ $v \in$
\{$1$, $ 2$, $ 3$, $ 5$, $ 6$, $ 7$, $ 10$, $ 11$, $ 15$, $ 21$, $ 22$, $ 23$, $ 26$, $ 29$, $ 35$, $ 38$, $ 39$, $ 42$, $ 70$, $ 71$, $ 74$, $ 86$, $ 95$, $ 134$, $ 182$, $ 186$, $ 230$, $ 266$, $ 326$, $ 470$, $ 506$, $ 806$, $ 854$\}.
\end{theorem}

\vspace{.5in}

We now consider the form $g:=(3,3,3,-2,2,2)$ of discriminant $\Delta =64$ and $|$Aut($g$)$|$=48. We see that $g(x,y,z):=3x^2+3y^2+3z^2 - 2yz+2xz+2xy$ is even if and only if $x+y+z \equiv 0 \imod{2}$. The substitution $x \mapsto (y+z),$  $y \mapsto (x-z),$ and $z \mapsto (x-y)$, ensures $x+y+z \equiv 0 \imod{2}$, and we have
\begin{equation}
\label{dis21}
\sum_{x+y+z \equiv 0 \imod{2}}q^{g(x,y,z)}=\sum_{x,y,z}q^{g(y+z,x-z,x-y)}=\sum_{x,y,z}q^{4(x^2+ y^2+z^2)}.
\end{equation}

The substitution $x \mapsto (y+z+1),$  $y \mapsto (x-z)$, and $z \mapsto (x-y)$, guarantees the condition $x+y+z \equiv 1 \imod{2}$, and we have
\begin{equation}
\label{dis211}
\sum_{x+y+z \equiv 1 \imod{2}}q^{g(x,y,z)}=\sum_{x,y,z}q^{g(y+z+1,x-z,x-y)}=\sum_{x,y,z}q^{(2x+1)^2 +(2y+1)^2 + (2z+1)^2}.
\end{equation}

Combining \eqref{dis21} and \eqref{dis211} yields the following lemma.
\begin{lemma}
\label{sif}
\begin{equation}
\sum_{x,y,z}q^{3x^2+3y^2+3z^2 - 2yz+2xz+2xy} = \sum_{x,y,z}q^{4(x^2+ y^2+z^2)} + \sum_{x,y,z}q^{(2x+1)^2 +(2y+1)^2 + (2z+1)^2}.
\end{equation}
\end{lemma}
Lemma \ref{sif} implies
\begin{equation}
\label{low3}
(3,3,3,-2,2,2;4n)=(1,1,1,0,0,0;n),
\end{equation}
\begin{equation}
\label{low3b}
(3,3,3,-2,2,2;8n+3)=(1,1,1,0,0,0;8n+3),
\end{equation}
and $(3,3,3,-2,2,2;n)=0$ for any $n \not\equiv 0,3,4 \imod{8}$.\\
A necessary condition that the integer $n$ is uniquely represented by $(3,3,3,-2,2,2)$ is 
\[
0<(3,3,3,-2,2,2;n) \leq |\mbox{Aut}((3,3,3,-2,2,2))|=48.
\]

We use $\text{Prelist}(1,1,1,0,0,0)$ to determine the solutions to 
\begin{equation}
\label{evenl2}
0<(3,3,3,-2,2,2;4n)=(1,1,1,0,0,0;n)\leq 48,
\end{equation}
and
\begin{equation}
\label{oddl2}
0<(3,3,3,-2,2,2;8n+3)=(1,1,1,0,0,0;8n+3)\leq 48.
\end{equation}

The integers $n$ which satisfy \eqref{evenl2} is exactly $\text{Prelist}(1,1,1,0,0,0)$, which we derived earlier in the section. The integers which satisfy \eqref{oddl2} is the finite subset of integers which are congruent to $3$ modulo $8$ and are in $\text{Prelist}(1,1,1,0,0,0)$.\\

Using Maple V.15 we check these solutions for unique representation, and find the following theorem.
\begin{theorem}
\label{333oe}
The form $(3,3,3,-2,2,2)$ uniquely represents the integer $n$ (up to action of automorphs) if and only if  $n=4^k\cdot v$, $k\geq 0, 4 \nmid v$, and \\ $v \in$
\{$3$, $ 4$, $ 8$, $ 11$, $ 19$, $ 20$, $ 24$, $ 35$, $ 40$, $ 43$, $ 52$, $ 56$, $ 67$, $ 84$, $ 88$, $ 91$, $ 115$, $ 120$, $ 148$, $ 163$, $ 168$, $ 184$, $ 232$, $ 235$, $ 280$, $ 312$, $ 372$, $ 403$, $ 427$, $ 532$, $ 568$, $ 760$, $ 1012$\}.
\end{theorem}

The three forms considered so far, $(1,1,1,0,0,0), (1,1,1,1,1,1),$ and $(3,3,3,-2,2,2)$, all have the maximum number of automorphs: 48. Let us briefly consider the form $(1,3,3,2,0,0)$ which has 8 automorphs.\\

The form $h=(1,3,3,2,0,0)$ is of discriminant $\Delta=32$ and $|$Aut($h$)$|$=8. To connect\\ $(1,3,3,2,0,0;n)$ to $(1,1,1,0,0,0;n)$ we note that $x^2 +3y^2 +3z^2 +2yz = x^2 +(y-z)^2 +2(y+z)^2$, and we can use Lemma \ref{kap1l} to reduce $(1,1,2,0,0,0)$ to $(1,1,1,0,0,0)$. Let $h(x,y,z)=x^2 +3y^2 +3z^2 +2yz$. We have

\begin{equation}
\label{133200}
\begin{aligned}
\sum_{x,y,z} q^{h(x,y,z)} &= \sum_{\substack{x,\\ y\equiv z \imod{2}}} q^{x^2 +y^2 +2z^2}\\
&=\sum_{x,y,z} q^{x^2 +2y^2 +8z^2}+\sum_{\substack{x,\\ y\equiv z \equiv 1 \imod{2}}} q^{x^2 +y^2 +2z^2}\\
&=\sum_{\substack{x \equiv 1\imod{2},\\y,z}} q^{x^2 +2y^2 +4z^2}+\sum_{x,y,z} q^{4(x^2 +y^2 +2z^2)}+\sum_{x \equiv y \equiv z \equiv 1 \imod{2}} q^{x^2 +y^2 +2z^2}\\
&=\sum_{\substack{x \equiv 1\imod{2},\\y,z}} q^{x^2 +2y^2 +4z^2}+\sum_{x,y,z} q^{8(x^2 +y^2 +z^2)}+\df{3}{2}\sum_{x \equiv y \equiv z \equiv 1 \imod{2}} q^{x^2 +y^2 +2z^2},\\
\end{aligned}
\end{equation}
where we employed the identity
\[
4\sum_{\substack{x \equiv 1\imod{2},\\y,z}} q^{4(x^2 +2y^2 +4z^2)} = \sum_{x \equiv y \equiv z \equiv 1 \imod{2}} q^{x^2 +y^2 +2z^2}.
\]
Making use of \eqref{kap3b} and observing that
\[
\sum_{x \equiv y \equiv z \equiv 1 \imod{2}}q^{x^2 +y^2 +2z^2} = 2\sum_{\substack{x\equiv y\equiv 1 \imod{2}\\z \equiv 0 \imod{2}}}q^{2(x^2 +y^2 +z^2)},
\]
we find
\begin{equation}
\label{l333}
(1,3,3,2,0,0;8n)=(1,1,1,0,0,0;n),
\end{equation}
\begin{equation}
\label{l333b}
(1,3,3,2,0,0;8n+4)=(1,1,1,0,0,0;4n+2),
\end{equation}
\begin{equation}
\label{l333c}
(1,3,3,2,0,0;4n+2)=0,
\end{equation}
\begin{equation}
\label{l333d}
6(1,3,3,2,0,0;2n+1)=(1,1,1,0,0,0;4n+2).
\end{equation}
Employing \eqref{l333} -- \eqref{l333d} along with Table \ref{prel} allows us to solve $0<(1,3,3,2,0,0;n) \leq 8$.\\

According to Table \ref{prel} we have $0<(1,1,1,0,0,0;n)\leq 8$ and $4\nmid n$ if and only if $n=1,3$. Hence \eqref{l333} implies $0<(1,3,3,2,0,0;8\cdot 4^k) \leq 8$ and $0<(1,3,3,2,0,0;24\cdot 4^k) \leq 8$. We see $(1,1,1,0,0,0;4n+2) \leq 48$ has the solutions $2n+1=1, 3, 5, 7, 9, 11, 15, 17, 21, 23, 29, 35, 39, 41, 51, 65, 71, 95$.\\

We now have all solutions to $0<(1,3,3,2,0,0;n) \leq 8$. Using Maple V.15 we directly check these integers for unique representation to arrive at the following theorem.

\begin{theorem}
\label{1333thm}
The form $(1,3,3,2,0,0)$ uniquely represents the integer $n$ (up to action of automorphs) if and only if  $n= 1,  3,  5,  7,  11,  15,  21,  23,  29,  35,  39,  71,$ or  $95$.
\end{theorem}


\section{Some ternary forms of discriminant $4096$ and $8192$}
\label{ex2}

The discriminant $4096$ is the largest discriminant which is a power of 4 and contains a ternary idoneal form \cite{jagykap}. Indeed the two forms $(5,13,20,-12,4,2)$ and $(5,12,20,8,4,4)$ are both idoneal and of discriminant $4096$. These forms are alike in the sense that they can be connected to $(1,1,1,0,0,0)$ to find the integers which they uniquely represent. The form $(5,12,20,8,4,4)$ has only 2 automorphs, and can be handled similarly to $(5,13,20,-12,4,2)$. We now deduce the integers uniquely represented by $(5,13,20,-12,4,2)$.\\

As mentioned, we can relate $(5,13,20,-12,4,2;n)$ to $(1,1,1,0,0,0;n)$ as given in the following theorem.
\begin{theorem}
\label{4096link}
Let $n$ be a nonnegative integer. We have
\begin{equation}
\label{zel}
(5,13,20,-12,4,2;64n)= (1,1,1,0,0,0;n),
\end{equation}
\begin{equation}
3(5,13,20,-12,4,2;32(2n+1))= (1,1,1,0,0,0;32(2n+1)),
\end{equation}
\begin{equation}
3(5,13,20,-12,4,2;16(4n+1))= (1,1,1,0,0,0;16(4n+1)),
\end{equation}
\begin{equation}
(5,13,20,-12,4,2;16(4n+3))= (1,1,1,0,0,0;16(4n+3)),
\end{equation}
\begin{equation}
3(5,13,20,-12,4,2;4(8n+5))= (1,1,1,0,0,0;8n+5),
\end{equation}
\begin{equation}
\label{endd}
12(5,13,20,-12,4,2;8n+5)= (1,1,1,0,0,0;8n+5),
\end{equation}
and $(5,13,20,-12,4,2;k)=0$ for any $k$ not covered by \eqref{4096link} -- \eqref{endd}.
\end{theorem}
Proof of the above theorem is elementary, but contains many details. Employing Theorem \ref{4096link} in conjunction with Table \ref{prel} gives all integers $n$ such that
 \[
(5,13,20,-12,4,2;n) \leq | \text{Aut}(5,13,20,-12,4,2)|=4,
\]
and so we can check which integers are uniquely represented by $(5,13,20,-12,4,2)$. In particular we point out that Table \ref{prel} implies there is no integer $n \equiv 0 \imod{64}$ with $(1,1,1,0,0,0,\tf{n}{64}) \leq 4$, and employing \eqref{zel} implies $(5,13,20,-12,4,2)$ does not uniquely represent infinitely many integers.\\
Instead of using Theorem \ref{4096link} to deduce the integers which are uniquely represented by \\$(5,13,20,-12,4,2)$, we employ the Siegel formula along with the method described in Section \ref{notprim}.\\

Using \eqref{mf} with $f:=(5,13,20,-12,4,2)$ we find
\begin{equation}
\label{thm4096}
(5,13,20,-12,4,2;n)=\dfrac{\sqrt{n}}{2\pi}\cdot d_{f,2}(n)\cdot L(1,\chi(n))\cdot P(n).
\end{equation}
Letting $n=4^av$, with $4\nmid v$, we see \eqref{thm4096} implies the bound
\begin{equation}
\label{thm22}
(5,13,20,-12,4,2;n)\geq \dfrac{2^{a-1}\sqrt{v}}{\pi}\cdot d_{f,2}(n)\cdot L(1,\chi(v)).
\end{equation} 
The local 2-adic density of $f$ is given by 
\begin{equation}
\label{4096ld1}
d_{f,2}(n) = \left\{ \begin{array}{ll}
				\frac{3}{2^{a-4}} & a \geq 3, ~v\equiv 1,2 \imod{4},\\
				\frac{1}{2^{a-5}} & a \geq 3, ~ v\equiv 3 \imod{8},\\
				0 &  v\equiv 7 \imod{8}.
     \end{array}
     \right.
\end{equation}
The values of $d_{f,2}(n)$ not covered in \eqref{4096ld1}, are listed in Table \ref{prelll}.
\begin{table}[htb]  
	\caption{} \label{prelll}
     \begin{center} 
  \begin{tabular}{ c | c | c | c  }
	&$a=0$& $a=1$& $a=2$ \\
	\hline
	   $v \equiv 1\imod{8}$  & 0& 0 & 4 \\ \hline
		 $v \equiv  3\imod{8}$ &0& 0 & 8  \\ \hline
		 $v \equiv 5\imod{8}$ & 4& 8 & 4  \\ \hline
     $v \equiv  2\imod{4}$ & 0& 0 & 4 . \\ 
    \hline
  \end{tabular}
\begin{flushright}
			 .
			 \end{flushright}
     \end{center}
   \end{table}

We used Sage 5.1 in computing the above densities. We refer the reader to \cite{berk} and \cite{hanke} for more details on computing local densities.\\

We can use \eqref{thm22} along with \eqref{4096ld1}, \eqref{h}, and Table \ref{prelll}, to find all $n$ that satisfy\\$0<(5,13,20,-12,4,2;n) \leq 4$. We write $n=4^av$ with $4\nmid v$ and split our analysis according to the cases $a=0,1;$ $a=2;$ or $a\geq 3$.\\

\textbf{\large Case 1:} $a=0,1.$\\

Since $a=0,1$ we see $(5,13,20,-12,4,2;n)=0$ unless $v \equiv 5 \imod{8}$. Hence we only consider $n=v$ or $n=4v$ with $v \equiv 5 \imod{8}$. In the case $n \equiv 5 \imod{8}$, we employ \eqref{h}, \eqref{thm22}, \eqref{4096ld1}, and Table \ref{prelll}, to find
\begin{equation}
\label{thm223}
(5,13,20,-12,4,2;n)\geq h(-4n).
\end{equation} 
We find $n=5, 13, 21, 37, 45, 85, 93, 133, 253$ are the solutions to $n \equiv 5 \imod{8}$ and $h(-4n) \leq 4$. Noting that $(5,13,20,-12,4,2;45)>4$, we see $n=45$ is the only spurious solution.

In the case $n=4v$ with $v \equiv 5 \imod{8}$, we employ \eqref{h}, \eqref{thm22}, \eqref{4096ld1}, and Table \ref{prelll}, to find
\begin{equation}
\label{thm224}
(5,13,20,-12,4,2;4v)\geq 4h(-4v).
\end{equation} 
There are no solutions to $h(-4v)=1$, and thus no solutions to $0<(5,13,20,-12,4,2;4v)\leq 4$ with $v \equiv 5 \imod{8}$.\\

\textbf{\large Case 2:} $a=2.$\\

Since $a=2$, we must consider $n=16v$ with $4 \nmid v$. Particularly, when $v=1$, we have $n=16$, and $(5,13,20,-12,4,2;16)=2$, so $n=16$ is a candidate for unique representation by $f$. However, $v=3$ implies $n=48$, and $(5,13,20,-12,4,2;48)>4$, so we need not consider $v=3$. Using \eqref{h}, \eqref{thm22}, \eqref{4096ld1}, and Table \ref{prelll}, we have
	\begin{equation}
	\label{p}
(5,13,20,-12,4,2;16v)\geq \left\{ \begin{array}{ll}
				24\cdot h(-v) &  3<v\equiv 3 \imod{8},\\
			4\cdot h(-4v) & 1< v\equiv 1,2 \imod{4}.\\
     \end{array}
     \right.
\end{equation}

Inequality \eqref{p} shows that we need only consider $v\equiv 1,2 \imod{4}$ when solving\\$0<(5,13,20,-12,4,2;16v) \leq 4$. We are left to solve $h(-4v)=1$ with $1<v\equiv 1,2 \imod{4}$, and we find the only solution is $v=2$. Therefore, we see that $n=32$ is a candidate for unique representation. \\

\textbf{\large Case 3:} $a\geq 3.$\\

We directly consider $n=4^a v$ with $4 \nmid v$ and $a \geq 3$. Using \eqref{h}, \eqref{thm22}, \eqref{4096ld1}, and Table \ref{prelll}, we have
\begin{equation}
\label{s}
(5,13,20,-12,4,2;4^a v)\geq \left\{ \begin{array}{ll}
6 &  v=1,\\
8 &  v=3,\\
				24h(-v) &  3<v\equiv 3 \imod{8},\\
			12h(-4v) & 1< v\equiv 1,2 \imod{4}.\\
     \end{array}
     \right.
\end{equation}

From \eqref{s}, we see there are no solutions to $(5,13,20,-12,4,2;4^a v) \leq 4$ and $a \geq 3$.\\

Combining the above case shows $\text{Prelist}(5,13,20,-12,4,2) =\{5, 13, 16, 21, 32, 37, 85, 93, 133, 253\}$. Employing Maple V.15 we easily check the elements of $\text{Prelist}(5,13,20,-12,4,2)$ for unique representation by $(5,13,20,-12,4,2)$ and we find the following theorem.

\begin{theorem}
\label{tt}
The form $(5,13,20,-12,4,2)$ uniquely represents $n$ (up to action of automorphs) if and only if $n=5, 13, 16, 21, 32, 37, 93, 133,$ or $253$.
\end{theorem}

We now treat the form $g:=(7,15,23,10,2,6)$. $g$ has discriminant $\Delta =2^{13} =8192$ and $|$Aut($g$)$|$=2. We remark that 8192 is the largest discriminant which is a power of two and contains a ternary idoneal form \cite{jagykap}. We can relate $(7,15,23,10,2,6;n)$ to $(1,1,1,0,0,0;n)$ as given in the following theorem.
\begin{theorem}
\label{8192link}
Let $n$ be a nonnegative integer. We have
\begin{equation}
\label{81}
(7,15,23,10,2,6;128n)=(1,1,1,0,0,0;n),
\end{equation}
\begin{equation}
3(7,15,23,10,2,6;64(2n+1))=(1,1,1,0,0,0;4n+2),
\end{equation}
\begin{equation}
3(7,15,23,10,2,6;32(4n+1))=(1,1,1,0,0,0;4n+1),
\end{equation}
\begin{equation}
(7,15,23,10,2,6;32(4n+3))=(1,1,1,0,0,0;4n+3),
\end{equation}
\begin{equation}
6(7,15,23,10,2,6;16(2n+1))=(1,1,1,0,0,0;4n+2),
\end{equation}
\begin{equation}
6(7,15,23,10,2,6;4(8n+7))=(1,1,1,0,0,0;2(8n+7)),
\end{equation}
\begin{equation}
\label{81b}
24(7,15,23,10,2,6;8n+7)=(1,1,1,0,0,0;2(8n+7)),
\end{equation}
and $(7,15,23,10,2,6;k)=0$ for any $k$ not covered by \eqref{81} -- \eqref{81b}.
\end{theorem}
 Employing Table \ref{prel} along with Theorem \ref{8192link} implies $(7,15,23,10,2,6)$ does not uniquely represent infinitely many integers. We chose to use the method of Section \ref{notprim} to treat the form $(7,15,23,10,2,6)$.

Using \eqref{mf} we have
\begin{equation}
\label{thm8192}
(7,15,23,10,2,6;n)=\dfrac{\sqrt{n}}{ 2\pi\sqrt{2}}d_{g,2}(n)\cdot L(1,\chi(2n))\cdot P(n,2).
\end{equation}

Let us write $n=4^av$, with $4\nmid v$. We find the 2-adic local densitiy of $g$ to be

\begin{equation}
\label{8192ld1}
d_{g,2}(n) = \left\{ \begin{array}{ll}
				\frac{3}{2^{a-5}} & a \geq 4, ~v\equiv 1 \imod{2},\\
				\frac{3}{2^{a-4}} & a \geq 4, ~ v\equiv 2 \imod{8},\\
				\frac{1}{2^{a-5}} & a \geq 4, ~ v\equiv 6 \imod{16},\\
				0 &  v\equiv 14 \imod{16}.
     \end{array}
     \right.
\end{equation}

The values of $d_{g,2}(n)$ not covered in \eqref{8192ld1} are listed in Table \ref{prellll}.

\begin{table}[htb]  
	\caption{} \label{prellll}
     \begin{center} 
  \begin{tabular}{ c | c | c | c|c  }
	&$a=0$& $a=1$& $a=2$& $a=3$\\
	\hline
	   $v \equiv 1,3,5\imod{8}$ & 0& 0 & 4 &4 \\ \hline
		$v \equiv  7\imod{8}$ & 4& 8 & 4 & 4 \\ \hline
		 $v \equiv 2\imod{8}$ & 0& 0 & 4 &6 \\ \hline
     $v \equiv  6\imod{16}$& 0& 0 & 8 &4. \\ 
    \hline
  \end{tabular}
\begin{flushright}
			 .
			 \end{flushright}
     \end{center}
   \end{table}

We now break our analysis into two cases depending on the parity of the order of 2 in $n$.\\

\textbf{\large Case 1:} $n=4^a \cdot v$ with $v$ odd.\\

Employing \eqref{h}, we have
\begin{equation}
\label{1l}
L(1,\chi(2n))=L(1,\chi(2v))=\frac{\pi}{2\sqrt{2v}}\cdot h(-8v),
\end{equation}
and combining \eqref{thm8192} with \eqref{1l} yields
\begin{equation}
\label{1thmm8192}
(7,15,23,10,2,6;n)\geq 2^{a-3}d_{g,2}(n)\cdot h(-8v).
\end{equation}

The form $(7,15,23,10,2,6)$ has 2 automorphs, so $0<(7,15,23,10,2,6;n)\leq 2$ is a necessary condition that $n$ be uniquely represented up to the action of automorphs.\\
Solving 
\begin{equation}
\label{m}
2^{a-3}d_{g,2}(n)\cdot h(-8v) \leq 2
\end{equation}
requires class number information up to class number 4. We remark that finding the solutions to \eqref{m} is similar to the process we used when we treated $(5,12,20,8,4,4)$ earlier in this section. We find $n=7, 15, 16, 23, 39, 71, 95$ are solutions to $0<(7,15,23,10,2,6;n)\leq 2$.\\

\textbf{\large Case 2:} $n=4^a \cdot 2\cdot v$ with $v$ odd.\\

In this case we see $L(1,\chi(2n))=L(1,\chi(v))$, and \eqref{thm8192} becomes
\begin{equation}
\label{thm88192}
(7,15,23,10,2,6;n)\geq\dfrac{\sqrt{n}}{ 2\sqrt{2}\pi}d_{g,2}(n)\cdot L(1,\chi(v)).
\end{equation}
Employing \eqref{h} and \eqref{thm88192} yields
\begin{equation}
\label{thmm81922}
(7,15,23,10,2,6;n)\geq \left\{ \begin{array}{ll}
        2^{a-3}\cdot d_{g,2}(n)&  v=1,\\
				&\\
				2^{a-2}\cdot d_{g,2}(n) &  v=3,\\
					&\\
					2^{a-2}\cdot d_{g,2}(n)\cdot h(-4v) &  1<v\equiv 1 \imod{4},\\
					&\\
				3\cdot 2^{a-2}\cdot d_{g,2}(n)\cdot h(-v) &  3<v\equiv 3 \imod{8},\\
					&\\
					2^{a-2}\cdot d_{g,2}(n)\cdot h(-v) &  v \equiv 7 \imod{8}.\\
     \end{array}
     \right.
\end{equation}

We employ \eqref{8192ld1}, \eqref{thmm81922}, Table \ref{prellll}, and the tables in the Appendix to find the only number $n=4^a \cdot 2\cdot v$ with $v$ odd, and $0<(7,15,23,10,2,6;n)\leq 2$, is $n=32$.

Combining Case 1 and Case 2 yields $\text{Prelist}(7,15,23,10,2,6) =\{7, 15, 16, 23, 32, 39, 71, 95\}$. Employing Maple V.15 we check the integers $7, 15, 16, 23, 32, 39, 71, 95,$ for unique representation and arrive at the following theorem.
\begin{theorem}
\label{t3t}
The form $(7,15,23,10,2,6) $ uniquely represents $n$ (up to action of automorphs) if and only if $n=7, 15, 16, 23, 32, 39, 71,$ or $95$.
\end{theorem}

\section{Resolving Some Conjectures of Kaplansky}
\label{ex3}

In the concluding remarks of \cite{kap}, Kaplansky regards $x^2 +y^2 +3z^2$ as ``the next challenge''. He computationally found the integers which are uniquely represented by $x^2 +y^2 +3z^2$. We now supply the proof of this conjecture.\\

In this section we deduce the integers which are uniquely represented by the forms $(1,3,3,0,0,0)$ and $(1,1,3,0,0,0)$. These two forms are intertwined with each other since for any nonnegative integer $n$, we have $(1,3,3,0,0,0;3n)=(1,1,3,0,0,0;n)$ and $(1,3,3,0,0,0;n)=(1,1,3,0,0,0;3n)$. Hence we also have $(1,3,3,0,0,0;n)=(1,3,3,0,0,0;9n)$ and $(1,1,3,0,0,0;n)=(1,1,3,0,0,0;9n)$.\\

Let $f:=(1,3,3,0,0,0)$ which is of discriminant 36 and has $|$Aut($f$)$|$=16. Equation \eqref{mf} gives
\begin{equation}
\label{thm36}
(1,3,3,0,0,0;n)=\dfrac{6}{\pi}\sqrt{n}\cdot d_{f,2}(n)\cdot d_{f,3}(n)\cdot L(1,\chi(9n))\cdot P(n,9),
\end{equation}
with all notation as defined in Section \ref{notprim}. Since $(1,3,3,0,0,0;n)=(1,3,3,0,0,0;9n)$, we only consider $9 \nmid n$. We write $n=4^a v$ with $4\nmid v$ and $9\nmid v$. We refer to \cite{berk} and \cite{hanke} for the following local density results:
\begin{equation}
\label{ld133}
d_{f,3}(n) = \left\{ \begin{array}{ll}
				2 & v\equiv 1\pmod{3}, \\
       0 & v\equiv 2\pmod{3}, \\
			\frac{4}{3} & v\equiv 3,6\pmod{9}, \\
     \end{array}
     \right.
\end{equation}
and
\begin{equation}
\label{lld2}
d_{f,2}(n) = \left\{ \begin{array}{ll}
				\frac{2^{a+2}-3}{2^{a+1}} & v\equiv 1,2\imod{4}, \\
				\vspace{-.2cm}&\\
       \frac{2^{a+1}-1}{2^{a}} & v\equiv 3\imod{8}, \\
			\vspace{-.2cm}&\\
			 2 &   v\equiv 7\imod{8}. \\
     \end{array}
     \right.
\end{equation}

Either by congruence considerations or by employing \eqref{ld133}, it is clear that $(1,3,3,0,0,0;n)=0$ when $n\equiv 2 \imod{3}$, so we do not consider such $n$. Using \eqref{el} we see
\begin{equation}
\label{l9}
L(1,\chi(9n))=\left\{ \begin{array}{ll}
\tf{4}{3} L(1,\chi(n)) & n\equiv 1 \imod{3},\\
\vspace{-.2cm}&\\
				L(1,\chi(n)) & n\equiv 3,6 \imod{9}.\\
     \end{array}
     \right.
\end{equation}
Employing \eqref{thm36}, \eqref{ld133}, and \eqref{l9}, we find
\begin{equation}
\label{n2}
(1,3,3,0,0,0;n)\geq \left\{ \begin{array}{ll}
\tf{16}{\pi}\sqrt{n}\cdot d_{f,2}(n)\cdot L(1,\chi(n)) & n\equiv 1 \imod{3},\\
\vspace{-.2cm}&\\
				\tf{8}{\pi}\sqrt{n}\cdot d_{f,2}(n)\cdot L(1,\chi(n)) & n\equiv 3,6 \imod{9}.\\
     \end{array}
     \right.
\end{equation}

Let $n=4^a v$ with $4\nmid v$, $9\nmid v$. Using \eqref{h} and \eqref{lld2} we have
\begin{equation}
\label{n3}
\tf{1}{\pi}\sqrt{n}\cdot d_{f,2}(n)\cdot L(1,\chi(n))= \left\{ \begin{array}{ll}
\tf{2^{a+2}-3}{8} & v=1,\\
\vspace{-.2cm}&\\
\tf{2^{a+1}-1}{2} & v=3,\\
\vspace{-.2cm}&\\
\tf{3(2^{a+1}-1)}{2}\cdot h(-v) & 3<v\equiv 3\imod{8},\\
\vspace{-.2cm}&\\
2^a \cdot h(-v) & v\equiv 7\imod{8},\\
\vspace{-.2cm}&\\
\tf{2^{a+2}-3}{4}\cdot h(-4v) & 1<v\equiv 1,2\imod{4}.\\
     \end{array}
     \right.
\end{equation}

Combining \eqref{n2} and \eqref{n3}, we find our lower bound for $(1,3,3,0,0,0;n)$ in terms of the class number, where $n=4^a \cdot v$ with $4\nmid v$, $9\nmid v$.
\begin{equation}
\label{n4}
(1,3,3,0,0,0;n)\geq \left\{ \begin{array}{ll}
2(2^{a+2}-3) & v=1,\\
\vspace{-.2cm}&\\
4(2^{a+1}-1) & v=3,\\
\vspace{-.2cm}&\\
24(2^{a+1}-1)\cdot h(-v) & v\equiv 19\imod{24},\\
\vspace{-.2cm}&\\
2^{a+4}\cdot h(-v) & v\equiv 7\imod{24},\\
\vspace{-.2cm}&\\
4(2^{a+2}-3)\cdot h(-4v) & 1<v\equiv 1,10\imod{12},\\
\vspace{-.2cm}&\\
12(2^{a+1}-1)\cdot h(-v) & 3<v\equiv 3,51\imod{72},\\
\vspace{-.2cm}&\\
2^{a+3}\cdot h(-v) & v\equiv 15,39\imod{72},\\
\vspace{-.2cm}&\\
2(2^{a+2}-3)\cdot h(-4v) & v\equiv 6,21,30,33 \imod{36}.\\
     \end{array}
     \right.
\end{equation}

The form $(1,3,3,0,0,0)$ has 16 automorphs. Employing \eqref{n4} along with the tables given in the Appendix, gives that there are exactly 53 numbers $n$ with $9 \nmid n$ and $0<(1,3,3,0,0,0;n)\leq 16$. These 53 numbers are the numbers in the set\\
$S=$\{1, 3, 4, 6, 7, 10, 12, 13, 15, 21, 22, 25, 30, 33, 34, 37, 42, 46, 57, 58, 66, 69, 70, 73, 78, 82, 85, 93, 97, 102, 105, 114, 130, 133, 138, 141, 142, 165, 177, 190, 193, 210, 213, 253, 258, 273, 282, 330, 345, 357, 438, 462, 498\}. We comment that we solving $0<(1,3,3,0,0,0;n)\leq 16$ requires class number information up to class number 8. We can use Maple V.15 to check the elements of $S$ for unique representation.
\begin{theorem}
\label{tt133}
The form $(1,3,3,0,0,0)$ uniquely represents $n$ (up to action of automorphs) if and only if $n=9^k\cdot v$, $k\geq 0$, with\\ $v\in$\{$1$, $ 3$, $ 6$, $ 10$, $ 13$, $ 21$, $ 22$, $ 30$, $ 33$, $ 34$, $ 37$, $ 42$, $ 46$, $ 57$, $ 58$, $ 66$, $ 69$, $ 78$, $ 82$, $ 85$, $ 93$, $ 102$, $ 114$, $ 130$, $ 138$, $ 141$, $ 142$, $ 165$, $ 177$, $ 190$, $ 210$, $ 213$, $ 253$, $ 258$, $ 282$, $ 345$, $ 357$, $ 462$, $ 498$\}.
\end{theorem}

From the comments at the beginning of this section, Theorem \ref{tt133} directly implies the following theorem.

\begin{theorem}
The form $(1,1,3,0,0,0)$ uniquely represents $n$ (up to action of automorphs) if and only if $n=9^k\cdot v$, $k\geq 0$, with\\ $v\in$\{$1$, $ 2$, $ 3$, $ 7$, $ 10$, $ 11$, $ 14$, $ 19$, $ 22$, $ 23$, $ 26$, $ 30$, $ 31$, $ 34$, $ 38$, $ 39$, $ 46$, $ 47$, $ 55$, $ 59$, $ 66$, $ 70$, $ 71$, $ 86$, $ 94$, $ 102$, $ 111$, $ 115$, $ 119$, $ 138$, $ 154$, $ 166$, $ 174$, $ 246$, $ 255$, $ 390$, $ 426$, $ 570$, $ 759$\}.
\end{theorem}

\noindent
We have now resolved the conjecture of Kaplansky, concerning the forms $(1,3,3,0,0,0)$ and\\ $(1,1,3,0,0,0)$. We comment that an almost identical analysis holds for the other idoneal forms of discriminant 36. We move on to treat a second conjecture of Kaplansky given in the concluding remarks of \cite{kap}.

In the concluding remarks of \cite{kap}, Kaplansky considers the forms $x^2+2y^2+3z^2$. His Theorem 7.2 states that the even integers which are uniquely represented by $x^2 +2y^2 +3z^2$ are the odd powers of 2. He does not show the proof of this theorem, but instead offers it as an exercise to the reader. Kaplansky computationally found the odd integers which are uniquely represented by $x^2 +2y^2 +3z^2$, but admits that the proof was not yet accessible.\\

We use the method of Section \ref{notprim} to find the odd integers which are uniquely represented by $x^2 +2y^2 +3z^2$. For the rest of our consideration of $(1,2,3,0,0,0;n)$ we take $n$ to be odd.\\

Let $g=(1,2,3,0,0,0)$ which is of discriminant 24 and has 8 automorphs. Employing \eqref{mf}, we find
\begin{equation}
\label{thm123}
(1,2,3,0,0,0;n)=\dfrac{18\sqrt{n}}{\pi\sqrt{6}}\cdot d_{g,2}(n)\cdot d_{g,3}(n)\cdot L(1,\chi(6n))\cdot P(n,6),
\end{equation}
with all notation as defined in Section \ref{notprim}.

\vspace{.3cm}

We find that for $n$ odd, we have $d_{g,2}(n)=1$. The 3-adic local density for $g$ is given by the following Lemma.
\begin{lemma}
\label{3lden}
Let $n=9^{b}v$ with $9\nmid v$. We have
\[
d_{g,3}(n) = \left\{ \begin{array}{ll}
				\frac{2(3^{b+1}-2)}{3^{b+1}} & v\equiv 1,2\imod{3}, \\
				\vspace{-.2cm}&\\
			2 & v\equiv 3\imod{9}, \\
			\vspace{-.2cm}&\\
				\frac{2(3^{b+1}-1)}{3^{b+1}} & v\equiv 6\imod{9}. \\
     \end{array}
     \right.
\]
\end{lemma}

Let $n=9^b\cdot v$ with $v$ odd and $9\nmid v$. Employing \eqref{h} we have
\begin{equation}
\label{l3}
L(1,\chi(6n))=\frac{\pi}{2\sqrt{6n}}\cdot h(-24n).
\end{equation}
Combining \eqref{3lden}, \eqref{thm123}, and \eqref{l3} we arrive at
\begin{equation}
\label{thmm4}
(1,2,3,0,0,0;n)\geq \left\{ \begin{array}{ll}
				\left(3-\tf{2}{3^b}\right)h(-24n) & v\equiv 1,2\imod{3}, \\
				\vspace{-.2cm}&\\
			3h(-24n) & v\equiv 3\imod{9}, \\
			\vspace{-.2cm}&\\
				\left(3-\tf{1}{3^b}\right)h(-24n) & v\equiv 6\imod{9}. \\
     \end{array}
     \right.
\end{equation}

The form $(1,2,3,0,0,0)$ has 8 automorphs, and so we solve for all odd $n$ with $0<(1,2,3,0,0,0;n)\leq 8$ which requires class number information up to 8. Utilizing \eqref{thmm4} along with the tables of the Appendix, we find the only odd $n$ with $0<(1,2,3,0,0,0;n)\leq 8$ to be
\[
n= 1, 3, 5, 7, 11, 13, 17, 19, 23, 35, 43, 47, 55, 73, 77, 83.
\]

	Maple V.15 can be used to check these 16 numbers for unique representation.
\begin{theorem}
\label{tt123}
The form $(1,2,3,0,0,0)$ uniquely represents odd $n$ (up to action of automorphs) if and only if $n=1,5,7,13,17,23,47,$ or $55$.
\end{theorem}

We have now confirmed and proven the observation of Kaplansky regarding $(1,2,3,0,0,0)$ \cite{kap}.

\section{Outlook}
\label{concl}

The method of employing the Siegel formula along with class number bounds easily extends to classifying the integers which are represented in essentially $k$ ways by an idoneal ternary quadratic form. One can use this paper as a guide to deduce the integers which are represented in an essentially unique way by any of the 794 idoneal ternary quadratic forms. It would be interesting to find the class number bounds that are necessary to classify the integers which are represented in an essentially unique way by any of the 794 idoneal ternary quadratic forms.\\

The restriction of the form being idoneal is not always necessary to find the integers which are uniquely represented by that form. To demonstrate this, we consider the form $(1,3,3,1,0,1)$ which is not idoneal, since $(1,1,11,1,1,1)$ shares the same genus. It can be shown that
\begin{equation}
\label{etap}
\vartheta(1,1,11,1,1,1,q) - \vartheta(1,3,3,1,0,1,q) = 4qE(q^4)^2E(q^{16}),
\end{equation}
where 
\[
E(q):= \prod_{n=1}^{\infty}(1-q^n).
\]

Equation \eqref{etap} shows that for $n\not\equiv 1 \imod{4}$ we have $(1,3,3,1,0,1;n)=(1,1,11,1,1,1;n)$. So when $n\not\equiv 1 \imod{4}$ we have
\begin{equation}
\label{cr1}
\df{(1,3,3,1,0,1;n)}{4} + \df{(1,1,11,1,1,1;n)}{12}= \df{(1,3,3,1,0,1;n)}{3},
\end{equation}
where we have used $|$Aut$(1,3,3,1,0,1)|=4$, and $|$Aut$(1,1,11,1,1,1)|=12$. Employing the Siegel formula gives\\
\begin{equation}
\label{a1}
	 (1,3,3,1,0,1; 32n) = (1,1,1,0,0,0;n),
	\end{equation}
	\begin{equation}
 (1,3,3,1,0,1; 16(2n+1) ) = (1,1,1,0,0,0;2(2n+1)),
	\end{equation}
\begin{equation}
 (1,3,3,1,0,1; 8(2n+1) ) = 0,
	\end{equation}
\begin{equation}
2(1,3,3,1,0,1; 4(2n+1) ) = (1,1,1,0,0,0;2(2n+1)),
	\end{equation}
\begin{equation}
(1,3,3,1,0,1; 2(2n+1) ) = 0,
	\end{equation}
\begin{equation}
\label{ao}
4(1,3,3,1,0,1; 4n+3 ) = (1,1,1,0,0,0;2(4n+3)).
	\end{equation}
   Using $|$Aut$(1,3,3,1,0,1)|=4$, Table \ref{prel}, and \eqref{a1} -- \eqref{ao}, we see $(1,3,3,1,0,1)$ does not uniquely represent any $n \not\equiv 1\imod{4}$. For $n \equiv 1 \imod{4}$, the Siegel formula gives
	\[
	3(1,3,3,1,0,1; n ) + (1,1,11,1,1,1;n) =  (1,1,1,0,0,0; 2n).
	\]
   Thus we obtain $0<(1,3,3,1,0,1; n) \leq  \tf{1}{3}(1,1,1,0,0,0; 2n)$ for $n \equiv 1 \imod{4}$. We are left to solve $(1,1,1,0,0,0; 2n) \leq 12$ for $n \equiv 1 \imod{4}$, and we find only $n=1$ as a solution. Indeed, $n=1$ is the only integer which is represented uniquely (up to the action of Aut$(1,3,3,1,0,1)$) by the form $(1,3,3,1,0,1)$.\\
   
    It is of interest to see which ternary quadratic forms uniquely represent (up to the action of automorphs) only a finite number of integers. The diagonal form $(1,2,4,0,0,0)$ considered by Kaplansky has this property, as well as the forms $(1,3,3,2,0,0), (5,13,20,-12,4,2), (7,15,23,10,2,6),$ and $(1,3,3,1,0,1)$, which were considered in this paper. The authors intend to address this topic in a subsequent paper.

      \section{Acknowledgements}
   \label{ack}

   We are grateful to William Jagy, Li--Chien Shen, John Voight, and Kenneth Williams for helpful discussions. We would like to thank George Andrews and James Sellers for their kind interest and encouragement. We are indebted to Keith Grizzell and Sue--Yen Patane for a careful reading of the manuscript and for many valuable suggestions.

\section{Appendix}

Below we list the negative of the 527 discriminants (fundamental and nonfundamental) of binary quadratic forms with class number $\leq 8$ according to the isomorphism class of the class group. We denote the class group of discriminant $D$ by $H(D)$, and we denote the cyclic group of order $n$ by $\mathbb{Z}_n$. We remark that we generated the tables below by utilizing the bounds found in \cite{watkins} along with equation \eqref{veq}. We then used PARI/GP V.2.7.0 to identify the isomorphism class of each class group.\\

\begin{table}[htb] 
\caption*{$H(D)\cong \mathbb{Z}_1 $}  \vspace{-.35cm}
\begin{tabular}{ccccccccccccc}
 3,& 4,& 7,& 8,& 11,& 12,& 16,& 19,& 27,& 28,& 43,& 67,& 163\\
\end{tabular}
   \end{table}

\begin{table}[htb] 
\caption*{$H(D)\cong \mathbb{Z}_2 $}  \vspace{-.35cm}
\begin{tabular}{cccccccccc}
 15,& 20,& 24,& 32,& 35,& 36,& 40,& 48,& 51,& 52,\\
 60,& 64,& 72,& 75,& 88,& 91,& 99,& 100,& 112,& 115,\\
 123,& 147,& 148,& 187,& 232,& 235,& 267,& 403,& 427&
\end{tabular}
   \end{table}
	
	\begin{table}[htb] 
\caption*{$H(D)\cong \mathbb{Z}_3 $}  \vspace{-.35cm}
\begin{tabular}{cccccccccc}
 23,& 31,& 44,& 59,& 76,& 83,& 92,& 107,& 108,& 124,\\
 139,& 172,& 211,& 243,& 268,& 283,& 307,& 331,& 379,& 499,\\
 547,& 643,& 652,& 883,& 907&&&&&
\end{tabular}
   \end{table}

	\begin{table}[htb] 
\caption*{$H(D)\cong \mathbb{Z}_2 \times \mathbb{Z}_2 $}  \vspace{-.35cm}
\begin{tabular}{cccccccccc}
 84,& 96,& 120,& 132,& 160,& 168,& 180,& 192,& 195,& 228,\\
 240,& 280,& 288,& 312,& 315,& 340,& 352,& 372,& 408,& 435,\\
 448,& 483,& 520,& 532,& 555,& 595,& 627,& 708,& 715,& 760,\\
 795,& 928,& 1012,& 1435&&&&&&
\end{tabular}
   \end{table}

	\begin{table}[htb] 
\caption*{$H(D)\cong \mathbb{Z}_4 $}  \vspace{-.35cm}
\begin{tabular}{cccccccccccccc}
 39,& 55,& 56,& 63,& 68,& 80,& 128,& 136,& 144,& 155,& 156,& 171,& 184,& 196,\\
 203,& 208,& 219,& 220,& 252,& 256,& 259,& 275,& 291,& 292,& 323,& 328,& 355,& 363,\\
 387,& 388,& 400,& 475,& 507,& 568,& 592,& 603,& 667,& 723,& 763,& 772,& 955,& 1003,\\
 1027,& 1227,& 1243,& 1387,& 1411,& 1467,& 1507,& 1555&&&&&&
\end{tabular}
   \end{table}

	\begin{table}[htb] 
\caption*{$H(D)\cong \mathbb{Z}_5 $}  \vspace{-.35cm}
\begin{tabular}{cccccccccc}
47,& 79,& 103,& 127,& 131,& 179,& 188,& 227,& 316,& 347,\\
412,& 443,& 508,& 523,& 571,& 619,& 683,& 691,& 739,& 787,\\
947,& 1051,& 1123,& 1723,& 1747,& 1867,& 2203,& 2347,& 2683&
\end{tabular}
   \end{table}

\begin{table}[htb] 
\caption*{$H(D)\cong \mathbb{Z}_6 $}  \vspace{-.35cm}
\begin{tabular}{cccccccccccccc}
87,& 104,& 116,& 135,& 140,& 152,& 175,& 176,& 200,& 204,& 207,& 212,& 216,\\
244,& 247,& 300,& 304,& 324,& 339,& 348,& 364,& 368,& 396,& 411,& 424,& 432,\\
436,& 451,& 459,& 460,& 472,& 484,& 492,& 496,& 515,& 531,& 540,& 588,& 628,\\
648,& 675,& 676,& 688,& 700,& 707,& 747,& 748,& 771,& 808,& 828,& 835,& 843,\\
856,& 867,& 891,& 931,& 940,& 963,& 988,& 1048,& 1059,& 1068,& 1072,& 1075,\\
1083,& 1099,& 1107,& 1108,& 1147,& 1192,& 1203,& 1219,& 1267,& 1315,& 1323,\\
1347,& 1363,& 1432,& 1563,& 1588,& 1603,& 1612,& 1675,& 1708,& 1843,& 1915,\\
1963,& 2227,& 2283,& 2403,& 2443,& 2515,& 2563,& 2608,& 2787,& 2923,& 3235,\\
3427,& 3523,& 3763,& 4075&&&&&&&
\end{tabular}
   \end{table}

	\begin{table}[htb] 
\caption*{$H(D)\cong \mathbb{Z}_7 $}  \vspace{-.35cm}
\begin{tabular}{cccccccccc}
71,& 151,& 223,& 251,& 284,& 343,& 463,& 467,& 487,& 587,\\
604,& 811,& 827,& 859,& 892,& 1163,& 1171,& 1372,& 1483,& 1523,\\
1627,& 1787,& 1852,& 1948,& 1987,& 2011,& 2083,& 2179,& 2251,& 2467,\\
2707,& 3019,& 3067,& 3187,& 3907,& 4603,& 5107,& 5923&&
\end{tabular}
   \end{table}

\begin{table}[htb] 
\caption*{$H(D)\cong \mathbb{Z}_8 $}  \vspace{-.35cm}
\begin{tabular}{ccccccccccccc}
95,& 111,& 164,& 183,& 248,& 272,& 295,& 299,& 371,& 376,& 380,& 392,& 395,\\
444,& 452,& 512,& 539,& 548,& 579,& 583,& 632,& 712,& 732,& 784,& 904,& 939,\\
979,& 995,&1024,& 1043,& 1156,& 1168,& 1180,& 1195,& 1252,& 1299,& 1339,& 1348,& 1528,\\
1552,& 1587,& 1651,& 1731,& 1795,& 1803,& 1828,& 1864,& 1912,& 1939,& 2059,& 2107,& 2248,\\
2307,& 2308,& 2323,& 2332,& 2395,& 2419,& 2587,& 2611,& 2827,& 2947,& 2995,& 3088,& 3283,\\
3403,& 3448,& 3595,& 3787,& 3883,& 3963,& 4195,& 4267,& 4387,& 4747,& 4843,& 4867,& 5587,\\
5707,& 5947,& 7987&&&&&&&&&&
\end{tabular}
 \end{table}

\begin{table}[htb] 
\caption*{$H(D)\cong \mathbb{Z}_4 \times \mathbb{Z}_2 $}  \vspace{-.35cm}
\begin{tabular}{ccccccccccccc}
224,& 260,& 264,& 276,& 308,& 320,& 336,& 360,& 384,& 456,& 468,& 504,& 528,\\
544,& 552,& 564,& 576,& 580,& 600,& 612,& 616,& 624,& 640,& 651,& 720,& 736,\\
768,& 792,& 819,& 820,& 832,& 852,& 868,& 880,& 900,& 912,& 915,& 952,& 987,\\
1008,& 1032,& 1035,& 1060,& 1128,& 1131,& 1152,& 1204,& 1240,& 1275,& 1288,& 1312,& 1332,\\
1360,& 1395,& 1408,& 1443,& 1488,& 1600,& 1635,& 1659,& 1672,& 1683,& 1752,& 1768,& 1771,\\
1780,& 1792,& 1827,& 1947,& 1992,& 2020,& 2035,& 2067,& 2088,& 2115,& 2128,& 2139,& 2163,\\
2212,& 2272,& 2275,& 2368,& 2392,& 2451,& 2475,& 2632,& 2667,& 2715,& 2755,& 2788,& 2832,\\
2907,& 2968,& 3172,& 3243,& 3355,& 3507,& 3627,& 3712,& 3843,& 4048,& 4123,& 4323,& 5083,\\
5467,& 6307&&&&&&&&&&&
\end{tabular}
 \end{table}

\begin{table}[htb] 
\caption*{$H(D)\cong \mathbb{Z}_2\times \mathbb{Z}_2\times \mathbb{Z}_2 $}  \vspace{-.35cm}
\begin{tabular}{cccccccccc}
420,& 480,& 660,& 672,& 840,& 960,& 1092,& 1120,& 1155,& 1248,\\
1320,& 1380,& 1428,& 1540,& 1632,& 1848,& 1995,& 2080,& 3003,& 3040,\\
3315&&&&&&&&&
\end{tabular}
\end{table}

\end{document}